\documentclass[12pt]{amsart}
\usepackage{amsfonts}
\usepackage{amsmath}
\usepackage{amsxtra}
\usepackage{amssymb,latexsym}
\usepackage[mathcal]{eucal}
\usepackage{amscd}
\input xy
\xyoption{all}

\newcommand{\Bcal}{\mathcal{B}}

\newcommand{\Hcal}{\mathcal{H}}

\newcommand{\Jcal}{\mathcal{J}}

\newcommand{\Lcal}{\mathcal{L}}

\newcommand{\Ncal}{\mathcal{N}}

\newcommand{\Rcal}{\mathcal{R}}

\newcommand{\Ucal}{\mathcal{U}}

\newcommand{\ch}{\mathbf{1}}

\newcommand{\N}{\mathbb{N}}

\newcommand{\al}{\alpha}

\newcommand{\ep}{\epsilon}

\newcommand{\tri}{\bigtriangleup}

\newcommand{\Aut}{{\rm{Aut\,}}}

\newcommand{\Homeo}{{\rm{Homeo}}}

\swapnumbers \theoremstyle{plain}
\newtheorem{thm}{Theorem}[section]
\newtheorem{cor}[thm]{Corollary}

\theoremstyle{definition}

\newtheorem{rmk}[thm]{Remark}

\numberwithin{equation}{section}

\begin{document}
\title[$\Aut(\mu)$ is Roelcke precompact]
{The group $\Aut(\mu)$ is Roelcke precompact}
\author{Eli Glasner}

\date{February 22, 2009}


\address{Department of Mathematics\\
     Tel Aviv University\\
         Ramat Aviv\\
         Israel}
\email{glasner@math.tau.ac.il}
\urladdr{http://www.math.tau.ac.il/~glasner/}




\thanks{Research partially supported by BSF (Binational USA-Israel)
grant no. 2006119.}

\begin{abstract}
Following a similar result of Uspenskij on the unitary group of 
a separable Hilbert space we show that with respect to the lower (or Roelcke)
uniform structure the Polish group $G= \Aut(\mu)$, of automorphisms
of an atomless standard Borel probability space $(X,\mu)$, is precompact.
We identify the corresponding compactification as the space of Markov
operators on $L_2(\mu)$ and deduce that the algebra of
right and left uniformly continuous functions, the algebra 
of weakly almost periodic functions, and the algebra of Hilbert 
functions on $G$, all coincide. 
Again following Uspenskij we also conclude that $G$ is totally minimal.  
\end{abstract}

\maketitle

Let $(X,\mu)$ be an atomless standard Borel probability space.
We denote by $\Aut(\mu)$ the Polish group of measure preserving
automorphisms of $(X,\mu)$ equipped with the weak topology.
If for $T \in G$ we let $U_T : L_2(\mu) \to L_2(\mu)$ be the 
corresponding unitary operator (defined by $U_Tf(x) = f(T^{-1}x)$),
then the map $T \mapsto U_T$ (the {\em Koopman map})
is a topological isomorphic embedding of
the topological group $G$ into the Polish topological group $\Ucal(H)$
of unitary operators on the Hilbert space $H = L_2(\mu)$ equipped with
the strong operator topology. The image of $G$ in $\Ucal(H)$
under the Koopman map is characterized as the 
collection of unitary operators $U \in \Ucal(H)$ for which $U(\ch)=\ch$ and
$Uf \ge 0$ whenever $f \ge 0$; see e.g. \cite[Theorem A.11]{G}.

It is well known (and not hard to see)   
that the strong and weak operator topologies coincide on $\Ucal(H)$
and that with respect to the weak operator topology, the group $\Ucal(H)$
is dense in the unit ball $\Theta$ of the space $\Bcal(H)$ of bounded 
linear operators on $H$.  
Now $\Theta$ is a compact space and as such it admits a unique uniform structure. The trace of the latter on $\Ucal(H)$ defines a uniform structure
on $\Ucal(H)$. We denote by $\Jcal$ the collection of Markov 
operators in $\Theta$, where $K \in \Theta$ is {\em Markov} if 
$K(\ch) = K^*(\ch) = \ch$ and $Kf \ge 0$ whenever $f \ge 0$.
It is easy to see that $\Jcal$ is a closed subset of $\Theta$.
Clearly the image of $G$ in $\Ucal(H)$ is contained in $\Jcal$ and it is 
well known that this image is actually dense in $\Jcal$ (see e.g. 
\cite{GK} or \cite{GLW}).
Thus, via the embedding of $G$ into $\Jcal$ we obtain also
a uniform structure on $G$. We will denote this uniform space
by $(G,\Jcal)$.

On every topological group $G$ there are two naturally defined
uniform structures $\Rcal(G)$ and $\Lcal(G)$. The {\em lower}
or  the {\em Roelcke} uniform structure on $G$ is defined as 
$\Ucal = \Rcal \wedge \Lcal$.
If $\Ncal$ is a base for the topology of $G$ at the neutral element
$e$, then with
$$
U_R= \{(x,y) : xy^{-1} \in U\} \qquad U_L=\{(x,y): x^{-1}y \in U\},
$$
the collections $\{U_R: U \in \Ncal\}$ and $\{U_L: U \in \Ncal\}$
constitute bases for $\Rcal(G)$ and $\Lcal(G)$ respectively.
A base for the Roelcke uniform structure is obtained by the
collection $\{V_U: U \in \Ncal\}$ where
$$
V_U = \{(x,y): \exists\ z,w \in U,\ y = zxw\}.
$$

In \cite{U} Uspenskij shows that the uniform structure induced from 
$\Theta$ on $\Ucal(H)$ coincides with the Roelcke structure of this
group. In this note we show that the same is true for $G = \Aut(\mu)$
and then, as in \cite{U},  deduce that $G$ is totally minimal
(see the definition below). 

The subject of Roelcke precompact (RPC) groups was
thoroughly studied by Uspenskij. 
In addition to the paper \cite{U} the interested reader can find 
more information about RPC groups in \cite{U01}, \cite{U02}
and \cite{U08}.
In \cite{U01} the author shows that the group $\Homeo(C)$ of 
self-homeomorphisms of the Cantor set $C$ with the compact-open 
topology is RPC. 
See also \cite[Sections 12 and 13]{GM} where the Polish
group $S_\infty(\N)$ of permutations of the natural numbers is shown 
to be RPC and where an alternative proof for the RPC property 
of $\Homeo(C)$ is indicated.

\section{$\Aut(\mu)$ is Roelcke precompact}

\begin{thm}\label{main}
The uniform structure induced from $\Jcal$ on $G$ coincides
with the Roelcke uniform structure $\Lcal \wedge \Rcal$.
Thus the Roelcke uniform structure on $G$ is precompact and the
natural embedding $G \to \Jcal$ is a realization of the Roelcke
compactification of $G$.
\end{thm}

\begin{proof}
Given $\ep >0$ and a finite measurable partition 
$\al=\{A_1,\dots,A_n\}$ of $X$ we set
$$
U_{\al,\ep}=\{ T \in G: \mu(A_i\tri T^{-1}A_i) < \ep, \ \forall\ 1 \le i \le n\},
$$ 
$$
W_{\al,\ep}= \{(S,T) \in G\times G: |\mu(A_i \cap S^{-1}A_j) 
- \mu(A_i \cap T^{-1}A_j) | < \ep,  \forall\ 1 \le i, j \le n\},
$$
and
$$
\tilde{W}_{\al,\ep}= \{(S,T) \in G\times G: \exists\
P, Q \in U_{\al,\ep},\ T = PSQ\}.
$$
Note that sets of the form $W_{\al,\ep}$ constitute a base
for the uniform structure on $G$ induced from $\Jcal$, while the
$\tilde{W}_{\al,\ep}$ form a base
for the Roelcke uniform structure on $G$.

If $T = PSQ$ with $P,Q \in U_{\al,\ep}$ then 
\begin{gather*}
|\mu(A_i \cap T^{-1}A_j) - \mu(A_i\cap S^{-1}A_j)|
= 
|\mu(A_i \cap Q^{-1}S^{-1}P^{-1}A_j) - \mu(A_i\cap S^{-1}A_j)| \\
 \le |\mu(A_i \cap Q^{-1}S^{-1}P^{-1}A_j) - \mu(Q^{-1}A_i \cap Q^{-1}S^{-1}P^{-1}A_j)|
 \\+
|\mu(A_i \cap S^{-1}P^{-1}A_j) - \mu(A_i \cap S^{-1}A_j)| < 2\ep.
\end{gather*}
Thus $\tilde{W}_{\al,\ep}\subset W_{\al,2\ep}$. This means that the
identity map $(G, \Lcal \wedge \Rcal) \to (G, \Jcal)$ is uniformly continuous.

For the other direction we start with a given $\tilde{W}_{\al,\ep}$.
Suppose $(S,T) \in W_{\al,\ep/{n^2}}$. Set
$$
A_{ij} = A_i \cap T^{-1}A_j, \quad A'_{ij} = A_i \cap S^{-1}A_j.
$$
We have $|\mu(A_{ij}) - \mu(A'_{ij})| < \ep/{n^2}$ for every $i$ and $j$.
Define a measure preserving $R \in G$ as follows. 
For each pair $i, j$ let $B_{ij} \subset A_{ij}$ with 
$\sum_{i,j} \mu(A_{ij} \setminus B_{ij}) < \ep$ and $\mu(B_{ij}) \le
\mu(A'_{ij})$. Next choose a measure preserving isomorphism $\phi_{ij}$
from $B_{ij}$ onto a subset $B'_{ij} \subset A'_{ij}$ and let
$\phi: \cup B_{ij} \to X$ be the map whose restriction to $B_{ij}$
is $\phi_{ij}$. Finally,  extend $\phi$ to an element $R \in G$
by defining $R$ on $X \setminus \cup B_{ij}$ to be any measure 
preserving isomorphism $X \setminus \cup B_{ij} \to X \setminus  \cup B'_{ij}$.

It is easy to check that $R$ is in $U_{\al,\ep}$, and for $TR^{-1}S^{-1}$ 
we have, up to sets of small measure,
\begin{align*}
TR^{-1}S^{-1}(A_i) &= 
TR^{-1}S^{-1}(\bigcup_{j=1}^n A_i \cap S A_j)\\
&=TR^{-1} (\bigcup_{j=1}^n S^{-1}A_i \cap  A_j)\\
& = T (\bigcup_{j=1}^n T^{-1}A_i \cap  A_j) = A_i.
\end{align*}
Thus also $TR^{-1}S^{-1}$ is in $U_{\al,\ep}$,
whence the equation $T = (TR^{-1}S^{-1}) S R$ implies 
$(S,T) \in \tilde{W}_{\al,\ep}$.
We have shown that  $W_{\al,\ep/{n^2}}\subset \tilde{W}_{\al,\ep}$ and it follows
that the identity map $(G,\Jcal) \to (G, \Lcal \wedge \Rcal)$ is also 
uniformly continuous.
\end{proof}

\begin{cor}
The Roelcke and the WAP compactifications of $G$ coincide.
Moreover, every bounded right and left uniformly continuous
function --- and hence also every WAP function --- on $G$
can be uniformly approximated by linear combinations of 
positive definite functions.
Or, in other words, every right and left uniformly continuous
function arises from a Hilbert representation.  
\end{cor}

\begin{proof}
Let $\Rcal o(G)$, denote the algebra of bounded right and left uniformly continuous complex-valued functions on $G$. We write $WAP(G)$ for the
algebra of weakly-almost-periodic complex-valued functions on $G$
and finally we let $\Hcal(G)$ be the algebra of Hilbert complex-valued functions 
on $G$ (i.e. the uniform closure of the algebra of all linear combinations of positive definite functions; the latter is also called the 
{\em Fourier-Stieltjes} algebra).
We then have
$$
\Rcal o(G) \supseteq WAP(G) \supseteq \Hcal(G).
$$
By Theorem \ref{main} these three algebras coincide for the topological group
$G=\Aut(\mu)$.
In fact, the functions of the form $F_{f}(T) = \langle Tf, f \rangle$ with
$f\in L_2(\mu)$ and $T \in \Jcal$, when restricted to $G$, 
are clearly positive definite. Since these functions generate the algebra 
$C(\Jcal)$, which by Theorem \ref{main} is canonically isomorphic to 
$\Rcal o(G)$, this shows that indeed $\Rcal o(G) = \Hcal(G)$.
\end{proof}

\begin{rmk}
By \cite{U} the same is true for the group $\Ucal(H)$.
For more details see \cite{M}. In the literature
a topological group $G$ for which
$$
WAP(G) = \Hcal(G)
$$
is called {\em Eberlein}. 
Thus both $\Ucal(H)$ and $\Aut(\mu)$ are Eberlein groups
and moreover for these groups
$$
\Rcal o(G) = WAP(G) = \Hcal(G).
$$
This fact for $\Ucal(H)$ was first shown in \cite{M}.
\end{rmk}

\section{$\Aut(\mu)$ is totally minimal}

A topological group is called {\em minimal} if it does not admit a strictly coarser Hausdorff group topology. It is {\em totally minimal} if all its 
Hausdorff quotient groups are minimal. 
Stoyanov proved that the unitary group $\Ucal(H)$ is totally 
minimal \cite{S}, \cite[Theorem 7.6.18]{DPS},  and Uspenskij provides in \cite{U}
an alternative proof based on his identification of $\Theta$ as the
Roelcke compactification of this group. Using Theorem \ref{main}
we have the following.

\begin{thm}\label{min}
The topological group $G=\Aut(\mu)$ is totally minimal. 
\end{thm}

For completeness we provide a proof of this theorem.
It follows Uspenskij's proof with some simplifications.
We will use though the next theorem of Uspenskij \cite[Theorem 3.2]{U}.

\begin{thm}\label{idem}
Let $S$ be a compact Hausdorff semitopological semigroup which satisfies 
the following assumption:
\begin{quote}
For every pair of idempotents $p, q \in S$ the conditions
$pq=p$ and $qp=p$ are equivalent.
(We write $p \le q$ when $p, q \in S$ satisfy these conditions.)
\end{quote}
Then every nonempty closed subsemigroup $K$ of $S$ contains
a least idempotent; i.e. an idempotent $p$ such that
$p \le q$ for every idempotent $q$ in $K$.
\end{thm}

It is not hard to check that $\Theta$ (and therefore also $\Jcal$) 
satisfies the assumption of this theorem.


\begin{proof}[A proof of Theorem \ref{min}]
Let $\tau$ denote the topology of a Hausdorff topological group
$G$ and suppose that $\tau'$ is a coarser Hausdorff group topology.
Then, the identity map $(G,\tau) \to (G,\tau')$ is continuous
and $\tau = \tau'$ iff this map is open. A moment's reflection now shows
that in order to prove that $G$ is totally minimal it suffices to check that 
every surjective homomorphism of Hausdorff topological groups 
$f: G \to G'$ is an open map.

So let $f: G \to G'$ be such a homomorphism and observe that
then $G'$ is Roelcke precompact and satisfies 
$\Rcal o(G') = WAP(G')$ as well. We denote by $\Jcal'$ the corresponding
(Roelcke and WAP) compactification of $G'$ and observe that
the dynamical systems $(\Jcal,G)$ and $(\Jcal',G')$ are their own
enveloping semigroups (see e.g. \cite{G}). 
Now, $(\Jcal,G)$ and $(\Jcal',G')$ being WAP systems,
the latter are compact semitopological semigroups. 
Moreover, the map $f: G \to G'$ naturally extends to a continuous 
homomorphism $F: \Jcal \to \Jcal'$. 
(This fact frees us from the need to use Proposition 2.1 from \cite{U}.)

Let $K = F^{-1}(e')$, where $e'$ is the neutral
element of $G'$. Clearly then $K$ is a closed subsemigroup of $\Jcal$.
Moreover, we have $gK =Kg = F^{-1}(g')$  whenever $g'=f(g) \in G'$.   
In fact, clearly $gK \subset F^{-1}(g')$ and if $F(q) = g'$ for some $q \in 
\Jcal$ then $F(g^{-1}q) = f(g^{-1})F(q) = g'^{-1}g'=e'$, hence
$g^{-1}q \in K$ and $q \in gK$. Thus $gK= F^{-1}(g')$ and symmetrically 
also $Kg= F^{-1}(g')$.  

Next observe that $gKg^{-1}=K$ for every $g \in G$. Thus if
$p$ is the least idempotent in $K$, provided by Theorem \ref{idem},
then $gpg^{-1}=p$ for all $g \in G$ and we conclude 
(an easy exercise) that either $p=I$,
the identity element of $G$, or $p$ is the projection on 
the space of constant functions (i.e. the operator of integration on 
$L_2(X,\mu)$).
In the second case we have $e'= F(p)=F(gp) = f(g)F(p) = f(g) e'= f(g)$
for every $g \in G$ and we conclude that $G'= \{e'\}$.

Suppose then that $p = I$. In that case we have $qK \subset K$
for every $q \in K$ and $qK$ being a closed subsemigroup, we
conclude that $I \in qK$. Similarly we get $I \in Kq$.
Whence $q$ is an invertible element of $\Jcal$, i.e. an element of $G$. 
Thus when $p=I$ we have $K \subset G$. 
 
Let now $g$ be an arbitrary element of $G$ and let $g' = f(g) \in G'$.
Suppose $G' \ni f(g_i) = g'_i \to g'$ 
is a convergent sequence in $G'$. With no loss in generality we
assume that $g_i$ converges to an element $q \in \Jcal$ and it
then follows that $F(q) = g'$. 
As $F^{-1}(g') = gK$ we conclude that $q \in gK \subset G$.
Now $f(g_i q^{-1}g) = f(g_i) = g'_i$ and $g_i q^{-1}g \to  g$.
This shows that $f$ is an open map and the proof is complete.
\end{proof}

\begin{rmk}
The group $G=\Aut(\mu)$ is in fact algebraically simple (Fathi \cite{F}). 
Thus minimality of $G$ implies total minimality. 
Note that with only slight changes the same proof applies to
$\Ucal(H)$ and thus we have here a simplified version of Uspenskij's
proof.
\end{rmk}

\begin{rmk}
It is perhaps worthwhile to mention here two other outstanding
properties of the Polish group $G=\Aut(\mu)$. The first, due to Giordano
and Pestov \cite{GP1} or \cite{GP2}, is that this group has the {\em 
fixed point on compacta property},
i.e. whenever $G$ acts on a compact space it admits a fixed point
(this property is also called {\em extreme amenability}).
The second is the fact that the natural unitary representation of $G$ on
$L_2(X,\mu)$ is irreducible (see \cite[Theorem 5.14]{G}).
\end{rmk}

{\it Acknowledgement:} Thanks are due to M. Megrelishvili, V. Pestov,
V. V. Uspenskij and B. Weiss who helped me to improve this note.
 
\bibliographystyle{amsplain}

\end{document}